\documentclass[11pt]{amsart}

\usepackage[T1]{fontenc}
\usepackage{lmodern}
\usepackage{microtype}
\usepackage{amsmath,amssymb,amsthm,mathtools}
\usepackage[new]{old-arrows}
\usepackage{enumitem}
\usepackage{hyperref}
\usepackage[nameinlink,capitalize,noabbrev]{cleveref}
\usepackage{comment}
\usepackage{mathrsfs}
\usepackage{lineno}
\hypersetup{
    colorlinks=true,
    linkcolor=blue,
    citecolor=blue,
    urlcolor=blue
}

\newtheorem{theorem}{Theorem}[section]
\newtheorem{proposition}[theorem]{Proposition}
\newtheorem{lemma}[theorem]{Lemma}
\newtheorem{corollary}[theorem]{Corollary}

\theoremstyle{definition}

\theoremstyle{remark}

\newcommand{\Z}{\mathbb Z}

\newcommand{\R}{\mathbb R}
\newcommand{\T}{\mathbb T}

\newcommand{\Isom}{\operatorname{Isom}}
\newcommand{\Fix}{\operatorname{Fix}}
\newcommand{\cd}{\operatorname{cd}}
\newcommand{\rank}{\operatorname{rank}}

\newcommand{\Ker}{\operatorname{Ker}}

\title[Group actions on aspherical Alexandrov spaces]
{Group Actions on Aspherical Alexandrov Spaces}

\author[N.~Bárcenas]{Noé Bárcenas}

\author[J.~Núñez-Zimbrón]{Jesús Núñez-Zimbrón}

\begin{document}

\begin{abstract}
Conner and Raymond proved that if a connected, compact Lie group $G$
acts effectively on a closed aspherical manifold $M$, then $G$ is a
torus and the fundamental group of $G$ injects into the center of
$\pi_1(M)$.
We prove the corresponding result for effective isometric actions on
closed, aspherical Alexandrov spaces.
\end{abstract}

\maketitle

% ============================================================
\section{Introduction}
% ============================================================

Let $M$ be a closed (that is, compact and without boundary), connected, \emph{aspherical} manifold, that is, a manifold such that its universal cover is contractible. Conner and Raymond proved in \cite{ConnerRaymond1970} (see also \cite[Section 3.1]{LeeRaymond}) that if a connected, compact Lie group $G$ acts effectively on $M$, then $G$ is a torus and the natural orbit homomorphism induces an injection $\pi_1(G)\hookrightarrow Z(\pi_1(M))$, where $Z(\pi_1(M))$ denotes the center of $\pi_1(M)$ (see also \cite{ConnerRaymond1971, ConnerRaymond1972} for subsequent work of Conner and Raymond in this direction).  Extensions of this result were obtained by Donnelly and Schultz for the so-called \emph{hyperaspherical} manifolds, that is, manifolds that admit degree-one maps onto aspherical manifolds \cite{DonnellySchultz1982}. Gottlieb, Lee and \"Ozaydin developed in \cite{GottliebLeeOzaydin1985} a broader theory of compact group actions and maps into \emph{$K$-manifolds} (closed, connected, oriented $n$-manifolds $M$ admitting a map $f\colon M\longrightarrow K(\pi,1)$ for some torsion-free group $\pi$, such that the induced homomorphism $f^*\colon H^n(\pi;\mathbb Z)\longrightarrow H^n(M;\mathbb Z)$ is surjective). Ku and Ku \cite{KuKu1983} and Tang \cite{Tang1993} studied related questions for \emph{aspherical $A_k(\pi)$-manifolds} (closed, connected, oriented manifolds $M$
admitting a map $f\colon M\to K(\pi,1)$ such that $f_*\colon H_k(M;\mathbb Q)\longrightarrow H_k(K(\pi,1);\mathbb Q)$ is nontrivial).

The purpose of this paper is to prove the following Alexandrov geometry analogue of the aforementioned theorem due to Conner and Raymond. Here, an Alexandrov space is a length space with a lower curvature bound in the sense of comparison triangles.  

\begin{theorem}
\label{thm:main}
Let $X$ be a closed, connected, aspherical Alexandrov $n$-space and $G$ a connected, compact Lie group acting effectively and isometrically on $X$. Then $G$ is a torus and there is an injective orbit homomorphism $\pi_1(G)\hookrightarrow Z(\pi_1(X))$. In particular, $\dim G \leq \rank_{\Z} Z(\pi_1(X))$.
\end{theorem}

We conjectured in \cite{BarcenasZimbron2021} that all closed, aspherical Alexandrov $3$-spaces are manifolds. Alattar and Freund \cite{AlattarFreund} confirmed this for spaces of non-negative curvature and arbitrary dimension $n\geq 3$. However, the conjecture turned out to be false in general: Kapovich \cite{Kapovich2025} constructed closed aspherical Alexandrov $n$-spaces which are not topological manifolds, also for $n\geq 3$. Thus, the class of compact aspherical Alexandrov spaces is strictly broader than that of aspherical manifolds and Theorem \ref{thm:main} becomes naturally relevant.

From Theorem \ref{thm:main}, when restricting to spaces with centerless fundamental group, one can drop the connectedness assumption on $G$ and obtain the following corollary.

\begin{corollary}
\label{cor:centerless}
Let $X$ be a closed, connected, aspherical Alexandrov space satisfying that $Z(\pi_1(X))=\{1\}$. Then every compact Lie group acting effectively and isometrically on $X$ is finite.
\end{corollary}

Our strategy to prove Theorem \ref{thm:main} has essentially two steps: First, we show that no nontrivial
finite $2$-group of isometries of $\widetilde X$ (the universal cover of $X$) can centralize the deck group $\pi_1(X)$ (cf. \cite[Lemma 3.1.13]{LeeRaymond}). The key observation here is that Smith theory and the geometry of fixed-point sets would yield an
upper bound for the cohomological dimension of $\pi_1(X)$
incompatible with the lower bound supplied by asphericity
and the top-dimensional mod-$2$ cohomology of $X$.
We then apply this obstruction to the lifted action of the universal covering group of $G$, which centralizes the deck group. It first rules out a nontrivial semisimple factor, forcing $G$ to be a torus, and then rules out a nontrivial kernel of the orbit homomorphism
$\pi_1(G)\to Z(\pi_1(X))$.  The proof is deliberately carried out over $\mathbb{Z}_2$ coefficients to avoid any orientability hypothesis.\\

\textbf{Statement on the use of AI}. OpenAI's Chat GPT-6 Astra was used to search the literature and audit our proofs for correctness. The authors take full responsibility for the mathematical content and the full manuscript.

\subsection{Acknowledgements}
The  first  named  author  thanks  the support  of  project PAPIIT  100326. Both authors acknowledge the support of project PAPIIT IA103925.
% ============================================================
\section{Preliminaries}
% ============================================================
In this section we collect some facts that we will need in the rest of the paper. For the sake of concision, we assume familiarity with the theory of Alexandrov spaces (of curvature bounded below) and refer the reader to the usual references \cite{AlexanderKapovitchPetrunin2024, BuragoBuragoIvanov2001}. We will however recall some facts here that are needed later. Throughout the article, all Alexandrov spaces will be assumed to be finite-dimensional unless otherwise stated. 

%================================
\subsection{Alexandrov spaces}
%================================

We begin by recalling that Alexandrov spaces are locally contractible on account of Perelman's conical neighborhood theorem \cite{Perelman1993}. In fact, Mitsuishi--Yamaguchi proved the stronger local Lipschitz
contractibility property \cite{MitsuishiYamaguchi2014}. It is also known that any Alexandrov space $X$ has the homotopy type of some CW-complex. Indeed as the topological dimension of $X$ coincides with its Hausdorff dimension (\cite[Corollary 6.5]{BuragoGromovPerelman1992}), a result of Borsuk (see for example \cite[Theorem 7.1]{Hu1965}) implies that $X$ is an ANR. Then, Milnor proved \cite[Theorem 1]{Milnor1959} that every ANR is dominated by a CW-complex. It is worth pointing out that it is not known whether Alexandrov spaces are actually homeomorphic to CW-complexes. 

Throughout the article we will deal with \emph{aspherical} spaces. A connected space $X$ is aspherical if its universal cover is contractible. For Alexandrov spaces this is equivalent to the condition that all higher homotopy groups vanish on account of them being dominated by CW-complexes. We prove this here in full for completeness sake. 

\begin{lemma}
\label{lem:aspherical}
Let $X$ be a connected Alexandrov space.  Then the following are
equivalent:
\begin{enumerate}[label=\textup{(\roman*)}]
    \item the universal cover $\widetilde X$ is contractible;
    \item $\pi_j(X)=0$ for every $j\geq2$.
\end{enumerate}
\end{lemma}

\begin{proof}
Let $p:\widetilde X\longrightarrow X$ be the universal covering map. Then the induced morphism $p_*:\pi_j(\widetilde X)      \longrightarrow \pi_j(X)$ is an isomorphism for $j\geq2$. If $\widetilde X$ is contractible, this immediately gives $\pi_j(X)=0$ for all $j\geq2$, proving the first implication. Conversely, suppose that $\pi_j(X)=0$ for all $j\geq2$.  Since $\widetilde X$ is simply connected, $\pi_1(\widetilde X)=0$ and the covering-space isomorphisms give that $\pi_j(\widetilde X)=0$ for $j\geq2$. Thus all higher homotopy groups of $\widetilde X$ vanish.

The lifted length metric makes $\widetilde X$ an Alexandrov space
with the same lower curvature bound than that of $X$. In particular,
$\widetilde X$ is locally contractible and has the homotopy type of a
CW complex. Therefore the map $\widetilde X\longrightarrow\{pt\}$
is a weak homotopy equivalence between spaces of CW type.  Whitehead's
theorem \cite[Theorem 4.5  in page 346]{Hatcher2002} then implies that it is a homotopy equivalence.  Hence $\widetilde X$ is contractible.
\end{proof}

Mitsuishi proved that a closed, connected Alexandrov $n$-space $X$ satisfies that $H^n(X,\mathbb{Z}_2)\cong \mathbb{Z}_2$ (see \cite[Lemma 2.6]{Mitsuishi2016}). We point out that Mitsuishi uses cohomology with compact supports, but as we  assume $X$ to be compact we can use singular cohomology. Furthermore, Mitsuishi proves this for the more general class of NB-spaces of which Alexandrov spaces are a particular case. 

A closed subset $E\subseteq X$ is \emph{extremal} if the following condition holds: Whenever $q\in X\setminus E$ and $p\in E$ is a local minimum of $d(q,\cdot)|_{E}$, the point $p$ is a critical point of $d(q,\cdot)$. That is, for every $v\in\Sigma_pX$, (where $\Sigma_pX$ denotes the space of directions at $p$) there exists a direction $\xi\in\Uparrow_p^q$ such that $\angle(v,\xi)\leq \frac{\pi}{2}$. Here, 
\[
\Uparrow_p^q
=
\left\{
\gamma^+(0):
\begin{array}{l}
\gamma\colon[0,d(p,q)]\to X \text{ is a unit-speed}\\
\text{minimizing geodesic from $p$ to $q$}
\end{array}
\right\}.
\]
We point out that a proper extremal subset of an Alexandrov $n$-space has dimension at most $n-1$, \cite[Lemma 5.3]{DengKapovitch}. Extremal subsets are locally contractible. Indeed, every open cover of an Alexandrov space admits a locally finite \emph{good refinement} \cite[Theorem~1.1(1)]{MitsuishiYamaguchi2019}. Here, a covering of a topological space is called \emph{good} if every nonempty finite intersection of elements in the covering is contractible. Moreover, Fujioka observes that, for every extremal subset $E$, the restriction $\{U_\alpha\cap E\}$ of a good covering $\{U_\alpha\}$ is a good cover of $E$ \cite[p.~5]{Fujioka2024}. It follows that every point of $E$ admits a basis of contractible relative neighborhoods.

If a compact Lie group $G$ acts isometrically on an Alexandrov space $X$, then the quotient $X/G$
is again an Alexandrov space with the same lower curvature bound, when equipped with the usual orbital metric. Furthermore, for every closed subgroup $H\leq G$, the image of $\Fix(H)$ in $X/G$ is an extremal subset (see \cite[\S4.2]{PerelmanPetrunin1993}; see also the explicit formulation in \cite[Example~2.3(3)]{Fujioka2025}). We now state a lemma we will need later.

\begin{lemma}
\label{lem:fixed-dimension}
Let $X$ be a connected Alexandrov $n$-space and let $F$ be a non-trivial finite group acting effectively and isometrically on $X$. If $\Fix(F)\neq \emptyset$, then $\dim \Fix(F)\leq n-1$. Moreover, $\Fix(F)$ is locally contractible.
\end{lemma}

\begin{proof}
Let $q:X\rightarrow X/F$ be the quotient map and set $E=q(\Fix(F))$. Then $E\subset X/F$ is an extremal subset. For $x,y\in \Fix(F)$, note that 
\[
d_{X/F}(q(x),q(y)) = \min_{g\in F}d_X(x,gy) =  d_X(x,y).
\]
Hence the restriction of the orbit map $q|_{\Fix(F)}:\Fix(F)\rightarrow E$ is distance-preserving.

We now claim that $E$ is a proper subset of $X/F$.  If we assume by contradiction that $E=X/F$ then every orbit would contain a point of $\Fix(F)$.  Given $x\in X$, choose $y\in \Fix(F)$ in the same orbit.  Then $x=gy$ for some $g\in F$, but $gy=y$ and thus $x=y$.  Hence $\Fix(F)=X$, contradicting the effectiveness of the action and the non-triviality of $F$. Therefore $E$ is a proper extremal subset and therefore $\dim E\leq n-1$. Finally, since $\Fix(F)$ and $E$ are isometric, $\dim \Fix(F)\leq n-1$. Local contractibility of $E$ follows from the local topological theory of extremal subsets described above, and the isometry $\Fix(F)\cong E$ transfers this property to $\Fix(F)$.
\end{proof}

% =======================================
\subsection{Cohomological dimension}
% =======================================

Let $k$ be a field and $\Gamma$ a discrete group. We let $\cd_k(\Gamma)$ be the projective dimension of the trivial $k[\Gamma]$-module $k$ and refer to it as the \emph{cohomological dimension of $\Gamma$ with coefficients in $k$}.
Equivalently,
\[
    \cd_k(\Gamma)
    =
    \sup\left\{
       q:
       H^q(\Gamma;V)\neq0
       \text{ for some }k[\Gamma]\text{-module }V
    \right\}.
\]
We refer to Brown \cite[Chapter~VIII]{Brown} for a review on the fundamental aspects of cohomological dimension.

We recall the form of the \emph{Cartan--Leray spectral sequence} that will be used below. Let $q:Y\longrightarrow Z$ be a connected regular covering with deck-transformation group $\Gamma$, and let $k$ be a field. Given a left $k[\Gamma]$-module $V$, let $\mathcal V$ denote the local coefficient system on $Z$ associated with
the covering $q$ and the $\Gamma$-module $V$. More explicitly, after choosing $z_0\in Z$ and $y_0\in q^{-1}(z_0)$, the covering determines a monodromy epimorphism
\[
\rho:\pi_1(Z,z_0)\longrightarrow\Gamma,
\]
and $\mathcal V$ is the local coefficient system determined by the
representation
\[
\pi_1(Z,z_0)
\xrightarrow{\rho}
\Gamma
\longrightarrow
\operatorname{Aut}_k(V).
\]
Its pullback $q^*\mathcal V$ to $Y$ is the constant local coefficient
system with fiber $V$.

The Cartan--Leray spectral sequence for the covering $q$, with these
coefficients, is
\[
E_2^{p,q}
=
H^p\!\left(
\Gamma;
H^q(Y;q^*\mathcal V)
\right)
\Longrightarrow
H^{p+q}(Z;\mathcal V).
\]
Here the $\Gamma$-action on
$H^q(Y;q^*\mathcal V)$ is induced by the deck-transformation action
together with the given action of $\Gamma$ on $V$.
See
\cite[Chapter~XVI, \S8, Theorem~8.4, p.~354;
\S9, p.~355, spectral sequence~(1)]{CartanEilenberg}.

% ======================================================
\subsection{Lifts of group actions}
\label{s:lifts}
%======================================================

 Let a connected Lie group $G$ act continuously on a connected, locally path connected, semilocally simply connected space $X$. Let $q:\widehat G\longrightarrow G$ be the universal cover group and $p:\widetilde X\longrightarrow X$ the universal cover of $X$. Then the action of $G$ on $X$ can be lifted to an action of $\widehat G$ on $\widetilde X$ as follows: We represent an element $\widehat g\in\widehat G$ by a path $g:[0,1]\longrightarrow G$, with $g(0)=e$. For $\widetilde x\in\widetilde X$, consider the path $\alpha(t)=g(t)\cdot p(\widetilde x)$ in $X$.  Let $\widetilde\alpha$ be its lift starting at $\widetilde x$. Then we define the action as $\widehat g\cdot\widetilde x := \widetilde\alpha(1)$. Homotopy lifting shows that this depends only on the class of the path $g(t)$ in $\widehat G$. Standard uniqueness of path lifting
shows that these transformations satisfy the group law. Hence they define an action of $\widehat G$ (cf. \cite[Theorem 9.1, Chapter I]{Bredon1972}). 

We now recall that the lifted action commutes with every deck transformation: Let $\gamma$ be a deck transformation. Then the path $\gamma\widetilde\alpha(t)$ is a lift of $\alpha(t)$ starting at $\gamma\widetilde x$. Uniqueness of path lifting therefore gives $\widehat g(\gamma\widetilde x)= \gamma(\widehat g\widetilde x)$. Thus $\widehat g\gamma=\gamma\widehat g$ for every $\widehat g\in\widehat G$ and every deck transformation $\gamma$.

Let us now impose more structure on $X$ and the action of $G$. Assume $X$ is a length space and we equip $\widetilde X$ with the lifted metric, making $p$ a local isometry. Then, it is easy to see that if the action of $G$ on $X$ is by isometries, then so is that of $\widehat G$ on $\widetilde X$.

To continue we now wish to recall the homomorphism $\pi_1(G)\to Z(\pi_1(X))$ used by Conner-Raymond \cite[Section 4]{ConnerRaymond1970}. For $\lambda\in\operatorname{ker}(q)$, the induced transformation of $\widetilde X$ covers the identity map of $X$, and hence is a deck transformation. We therefore obtain a homomorphism
\[
\rho: \operatorname{ker}(q) \longrightarrow \operatorname{Deck}(\widetilde X/X),
\]
which under the standard identifications $\operatorname{ker}(q)\cong\pi_1(G)$ and $\operatorname{Deck}(\widetilde X/X)\cong\pi_1(X)$ gives us a homomorphism 
\[
\rho: \pi_1(G)\longrightarrow\pi_1(X).
\]
The homomorphism $\rho$ agrees, up to the convention used to identify deck transformations with elements of the fundamental group, with the
homomorphism
\[
(\operatorname{ev}_x)_*\colon
\pi_1(G)\longrightarrow\pi_1(X)
\]
induced by the orbit map evaluation 
\[
\operatorname{ev}_x\colon G\longrightarrow X,
\qquad
g\longmapsto g\cdot x,
\]
where $x\in X$ is any chosen basepoint. Finally, since the lifted action of $\widehat G$ on $\widetilde X$ commutes with every deck transformation, we have that $\rho(\pi_1(G))\subseteq Z(\pi_1(X))$. 

% ======================================================
\subsection{Smith theory}
% ======================================================
Smith theory studies the relationship between the topology of a space and that of the fixed-point set of an action by a finite $p$-group. Here, a finite group is a \emph{$p$-group} if its order is a power of $p$. In particular, Smith theory provides restrictions on the mod-$p$ (\v Cech) homology and cohomology of the fixed-point set in terms of those of the ambient space. We refer the reader to Bredon's book \cite[Chapter VII]{Bredon1972} for an exposition of the subject. 

A paracompact Hausdorff space is called \emph{finitistic} if it satisfies the so called \emph{Swan condition}, that is, if every open cover admits a finite-dimensional
refinement, and it is clear that every finite-dimensional space is finitistic (see
\cite[Chapter III, \S6]{Bredon1972}).

\begin{theorem}
\label{thm:smith}
Let $F$ be a finite $p$-group acting continuously on a finitistic space $Y$. If $Y$ is mod-$p$ acyclic in \v{C}ech cohomology, that is, 
\[
\check H^i(Y;\mathbb{Z}_p)
\cong
\begin{cases}
\mathbb{Z}_p, & i=0,\\
0, & i>0,
\end{cases}
\]
then $\Fix(F)$ is nonempty and mod-$p$ acyclic in \v{C}ech cohomology:
\[
\check H^i(\Fix(F);\mathbb Z_p)
\cong
\begin{cases}
\mathbb Z_p, & i=0,\\
0, & i>0.
\end{cases}
\]
\end{theorem}

This is the case $n=0$ of \cite[Chapter III, \S7, Theorem 7.11]{Bredon1972}.
Indeed, a mod-$p$ acyclic space is a mod-$p$ \v{C}ech cohomology $0$-disk. The cited theorem asserts that, if $F$ is a finite $p$-group and $Y$ is finitistic, then the fixed-point set of a
mod-$p$ \v{C}ech cohomology $n$-disk is a mod-$p$ \v{C}ech cohomology $r$-disk for some $0\leq r\leq n$. When $n=0$, necessarily $r=0$, which gives the assertion.

The preceding theorem applies, in particular, to finite $p$-group
actions by isometries on finite-dimensional Alexandrov spaces as they are finitistic.

Suppose now that $X$ is a contractible finite-dimensional Alexandrov space and that a finite $p$-group $F$ acts isometrically on $X$. Since $X$ is contractible, 
\[
\check H^i(X;\mathbb Z_p)
\cong
\begin{cases}
\mathbb Z_p, & i=0,\\
0, & i>0.
\end{cases}
\]
Hence Theorem~\ref{thm:smith} implies that $\Fix(F)$ is nonempty and mod-$p$ acyclic in \v{C}ech cohomology.

In the applications below we shall use this conclusion in terms
of ordinary singular homology. First, $\operatorname{Fix}(F)$
is locally contractible. Indeed, let $F_{\mathrm{eff}}$ be the
image of $F$ in $\operatorname{Isom}(X)$. If $F_{\mathrm{eff}}$ is trivial, then $\operatorname{Fix}(F)=X$,
which is locally contractible. Otherwise, $F_{\mathrm{eff}}$
is a nontrivial finite group acting effectively and isometrically
on $X$, and $\operatorname{Fix}(F)=\operatorname{Fix}(F_{\mathrm{eff}})$. Since this fixed-point set is nonempty by Smith theory,
Lemma \ref{lem:fixed-dimension} implies that it is locally contractible.
Since $\operatorname{Fix}(F)$ is metrizable and locally contractible, its \v{C}ech and singular cohomology groups agree by \cite[Corollaries~6.8.8 and~6.9.5]{Spanier}.

 Consequently,
\[
H^i(\Fix(F);\mathbb Z_p)
\cong
\begin{cases}
\mathbb Z_p, & i=0,\\
0, & i>0.
\end{cases}
\]

Finally, since $\mathbb Z_p$ is a field, the universal coefficient theorem gives
\[
H^i(\Fix(F);\mathbb Z_p)
\cong
\operatorname{Hom}_{\mathbb Z_p}
\bigl(H_i(\Fix(F);\mathbb Z_p),\mathbb Z_p\bigr).
\]
It follows that $H_i(\Fix(F);\mathbb Z_p)=0$ for every $i>0$. Looking now at $i=0$, we have that  $H^0(\Fix(F);\mathbb Z_p)\cong\mathbb Z_p$ implies that $\Fix(F)$ has exactly one path-component. Hence $H_0(\Fix(F);\mathbb Z_p)\cong\mathbb Z_p$.  Thus $\Fix(F)$ is mod-$p$ acyclic also with respect to singular homology.

% ============================================================
\section{Cohomological dimension estimates}
% ============================================================

Let $X$ be a closed, connected, aspherical Alexandrov $n$-space and let $\pi=\pi_1(X)$. In this section we will first produce a lower estimate of the cohomological dimension of $\pi$ with coefficients in $\mathbb{Z}_2$. 

% ------------------------------------------------------------
\subsection{The lower bound}
% ------------------------------------------------------------

We first deal with the lower bound on $\cd_{\mathbb{Z}_2}(\pi)$. We formulate the following result for Alexandrov spaces but point out that we only really require that the space admits a mod-$2$ fundamental class (on top of the other topological assumptions). 

\begin{lemma}
\label{lem:cd-lower}
Let $X$ be a closed, connected, aspherical Alexandrov $n$-space and let $\pi=\pi_1(X)$. Then $\cd_{\mathbb{Z}_2}(\pi)\geq n$.
\end{lemma}

\begin{proof}
Since $X$ is aspherical, it is an Eilenberg-Maclane space $K(\pi,1)$. Since $H^n(X;\mathbb Z_2)\cong\mathbb Z_2$, the universal
coefficient theorem over $\mathbb Z_2$ implies that
$H_n(X;\mathbb Z_2)\neq 0$. Therefore, $H_n(\pi;\mathbb{Z}_2) \cong H_n(X;\mathbb{Z}_2) \cong \mathbb{Z}_2 \neq0$. Let us assume by contradiction that  $\cd_{\mathbb{Z}_2}(\pi)<n$. Then the trivial $\mathbb{Z}_2[\pi]$-module
$\mathbb{Z}_2$ would admit a projective resolution of length less than $n$. Consequently
\[
    \operatorname{Tor}^{\mathbb{Z}_2[\pi]}_n(\mathbb{Z}_2,\mathbb{Z}_2)=0.
\]
However, $\operatorname{Tor}^{\mathbb{Z}_2[\pi]}_n(\mathbb{Z}_2,\mathbb{Z}_2)\cong H_n(\pi;\mathbb{Z}_2)$, which is a contradiction.  
\end{proof}

% ------------------------------------------------------------
\subsection{The upper bound}
% ------------------------------------------------------------
We are now interested in producing an upper bound on the cohomological dimension of deck transformation groups of covers of (subsets of) Alexandrov spaces. There is a standard way to obtain an upper bound for the cohomological dimension of a group from a free action on an acyclic CW complex. Namely, if $Y$ is a $d$-dimensional, $k$-acyclic free
$\Gamma$-CW complex, then its augmented cellular chain complex
\[
0 \longrightarrow C_d(Y;k) \longrightarrow \cdots
\longrightarrow C_0(Y;k) \longrightarrow k
\longrightarrow 0
\]
is a free $k[\Gamma]$-resolution of the trivial module $k$.
Indeed, each $C_i(Y;k)$ is a free $k[\Gamma]$-module,
with one generator for each $\Gamma$-orbit of $i$-cells,
and $k$-acyclicity makes the augmented complex exact.
Hence $\operatorname{cd}_k(\Gamma)\leq d$.

We cannot apply this argument directly in our setting. The space $Y$
that will arise below is isometric to an extremal subset of an Alexandrov space, and the estimate available to us is an estimate of its covering dimension. We do not have a $\Gamma$-equivariant CW structure on $Y$ (or even a CW structure at all) whose cellular dimension is controlled by this
covering dimension. The following lemma replaces the cellular
argument by the Cartan--Leray spectral sequence and the
cohomological-dimension theorem for finite-dimensional paracompact
spaces.

\begin{lemma}
\label{lem:CL-dimension}
Let $k$ be a field, and let $q:Y\longrightarrow Z$ be a connected regular covering with deck-transformation group $\Gamma$. Assume that

\begin{enumerate}
    \item $Y$ is $k$-acyclic,
    \item $Z$ is paracompact, Hausdorff, and locally contractible,
    \item the covering dimension of $Z$ satisfies $\dim Z\leq d<\infty$.
\end{enumerate}
Then $\operatorname{cd}_k(\Gamma)\leq d$.
\end{lemma}

\begin{proof}
Let $V$ be an arbitrary left $k[\Gamma]$-module. It suffices to prove $H^p(\Gamma;V)=0$ for every $p>d$. Let $\mathcal V$ be the local coefficient system on $Z$ associated with the regular covering $q:Y\to Z$ and the $\Gamma$-module $V$. The Cartan--Leray spectral sequence is then 
\[
    E_2^{p,q}
      =
    H^p\!\left(\Gamma;
       H^q(Y;q^*\mathcal V)\right)
      \Longrightarrow
    H^{p+q}(Z;\mathcal V).
\]
The pullback local system $q^*\mathcal V$ is constant with fiber $V$,
since the monodromy representation becomes trivial after restriction to
$q_{\ast}\pi_1(Y,y_0)=\ker\rho$. Hence
\[
    H^q(Y;q^*\mathcal V)\cong H^q(Y;V),
\]
where $V$ on the right denotes the constant coefficient system.

Since $V$ is a $k$-vector space and $k$ is a field, the universal
coefficient theorem gives
\[
    H^q(Y;V)
       \cong
    \operatorname{Hom}_k(H_q(Y;k),V).
\]
Since $Y$ is $k$-acyclic, $H^q(Y;V)=0$ for $q>0$, while $H^0(Y;V)\cong V$. Under this identification, the $\Gamma$-action on $H^0(Y;V)$ appearing
in the Cartan--Leray spectral sequence is the original
$k[\Gamma]$-module structure on $V$. Thus
\[
    E_2^{p,0}=H^p(\Gamma;V),
    \qquad
    E_2^{p,q}=0\quad(q>0).
\]
The spectral sequence therefore collapses, and for every $p\geq0$ we
obtain a natural isomorphism
\begin{equation}
\label{eq:CL-collapse}
    H^p(\Gamma;V)
       \cong
    H^p(Z;\mathcal V).
\end{equation}

Now, let $\mathscr{V}$ denote the locally constant sheaf of $k$-modules associated
with the local system $\mathcal V$. Since $Z$ is locally contractible,
it is homologically locally connected. Since it is also paracompact and
Hausdorff, the comparison theorem between singular cohomology with local
coefficients and sheaf cohomology gives natural isomorphisms
\[
    H^p(Z;\mathcal V)
       \cong
    H^p_{\mathrm{sh}}(Z;\mathscr{ V});
\]
see
\cite[Chapter III, \S1, Theorem~1.1]{BredonSheaf}.

On the other hand, for a paracompact Hausdorff space the sheaf
cohomological dimension over $k$ is bounded above by its
covering dimension. More precisely,
\cite[Chapter II, \S16, Corollary~16.34]{BredonSheaf}
implies, from $\dim Z\leq d$, that
\[
    H^p_{\rm sh}(Z;\mathscr F)=0
    \qquad (p>d)
\]
for every sheaf $\mathscr {F}$ of $k$-modules on $Z$. In particular this is the case for $F=V$. Consequently,
\[
    H^p(Z;\mathcal V)=0
    \qquad(p>d),
\]
and \eqref{eq:CL-collapse} yields
\[
    H^p(\Gamma;V)=0
    \qquad(p>d).
\]

Since $V$ was an arbitrary left $k[\Gamma]$-module, it follows that $\operatorname{cd}_k(\Gamma)\leq d$ as we wanted. 
\end{proof}

% ============================================================
\section{Proofs of Theorem \ref{thm:main} and Corollary \ref{cor:centerless}}
% ============================================================

To prove Theorem \ref{thm:main} we will first use Smith theory to prove the following result. In what follows we denote the fundamental group $\pi_1(X)$ of a connected space $X$ by $\pi$. 

\begin{proposition}
\label{prop:centralizer}
Let $X$ be a closed, connected, aspherical Alexandrov $n$-space and $p:\widetilde X\rightarrow X$ its universal cover.  Then there is no non-trivial finite $2$-subgroup $F\leq\Isom(\widetilde X)$ which centralizes $\pi$, that is, such that its elements commute with every deck transformation.
\end{proposition}

\begin{proof}
Let us suppose to the contrary that such a group $F$ exists. Set $Y:=\Fix(F)$.  Since $X$ is aspherical, $\widetilde X$ is contractible.  In particular it is mod-$2$ acyclic. Smith theory (Theorem \ref{thm:smith}) therefore implies that $Y$ is nonempty and mod-$2$ acyclic, which in terms of reduced homology says that $\widetilde H_*(Y;\mathbb{Z}_2)=0$.

By Lemma \ref{lem:fixed-dimension}, $\dim Y\leq n-1$, and $Y$ is locally contractible. Because $F$ centralizes $\pi$, $Y$ is $\pi$-invariant. Indeed, if $y\in Y$, $\gamma\in\pi$, and $a\in F$, then $a(\gamma y)=\gamma(ay)=\gamma y$. Thus $\gamma y\in Y$. The restricted action $\pi\curvearrowright Y$ is free and properly discontinuous since it is the restriction of the deck transformations action on $\widetilde X$.

Set $Z:=p(Y)$, equipped with the subspace topology inherited from $X$. We first observe that $p^{-1}(Z)=Y$. Indeed, if $p(z)=p(y)$ for some $y\in Y$, then $z=\gamma y$ for some $\gamma\in\pi$, and the $\pi$-invariance of $Y$ implies that $z\in Y$.

The restriction $q:=p|_Y:Y\to Z$ is a covering map.
To see this, let $U\subseteq X$ be an evenly covered open
neighborhood, with
\[
p^{-1}(U)=\coprod_\alpha U_\alpha.
\]
The equality $p^{-1}(Z)=Y$ implies that each homeomorphism $p|_{U_\alpha}:U_\alpha\to U$ restricts to a homeomorphism $U_\alpha\cap Y \longrightarrow U\cap Z$.  Thus $U\cap Z$ is evenly covered by $q$.
In particular, $q$ is an open quotient map. Its fibers are
precisely the $\pi$-orbits, so it induces a homeomorphism
$Y/\pi\cong Z$.

Since $Y$ is mod-$2$ acyclic in singular homology, it is
path connected. The restricted action of $\pi$ on $Y$ is
effective: an element acting trivially on the nonempty set
$Y$ fixes a point of $\widetilde X$, and hence is the identity
deck transformation. Moreover, $\pi$ acts transitively on
each fiber of $q$. Every deck transformation of $q$ therefore
agrees at one point with the restriction of some
$\gamma\in\pi$. Uniqueness of lifts on the connected space
$Y$ implies that they agree everywhere. Consequently, $q$
is a connected regular covering with deck group $\pi$.

The set $Y$ is closed in $\widetilde X$, being an intersection
of fixed-point sets of isometries. Since $p^{-1}(Z)=Y$ and
$p$ is a quotient map, $Z$ is closed in $X$. Thus $Z$ is
compact metrizable, and in particular paracompact Hausdorff.
The local homeomorphisms supplied by $q$ also show that $Z$
is locally contractible, because $Y$ is locally contractible.

Finally, the proof of Lemma \ref{lem:fixed-dimension} identifies $Y$ isometrically with a proper extremal subset of $\widetilde X/F$.
Hausdorff and covering dimensions agree for extremal subsets
(see \cite{Fujioka2025}), so Lemma \ref{lem:fixed-dimension} also gives
$\dim_{\mathrm H}Y\le n-1$. 
Since $p$ is $1$-Lipschitz for the lifted length metric,
\[
\dim Z
\le \dim_{\mathrm H}Z
\le \dim_{\mathrm H}Y
\le n-1.
\]
All hypotheses of Lemma \ref{lem:CL-dimension} therefore hold, with
$k=\mathbb Z_2$, and give $\operatorname{cd}_{\mathbb Z_2}(\pi)\le n-1$.
This contradicts Lemma \ref{lem:cd-lower}, which gives $\operatorname{cd}_{\mathbb Z_2}(\pi)\ge n$. Therefore no such nontrivial finite $2$-subgroup $F$ exists.

\end{proof}

With this we are ready to prove Theorem \ref{thm:main}.

\begin{proof}[Proof of Theorem \ref{thm:main}]
Let $q:\widehat G\rightarrow G$ be the universal covering group of $G$. By the structure theory of compact connected Lie groups (see for example \cite[Theorem 11.2]{SanMartin}),  $\widehat G\cong \R^k\times S$, where $S$ is a compact, connected, simply-connected semisimple Lie group.  Moreover,
$\Lambda:=\Ker q$ is a discrete central subgroup of $\widehat G$,  $\widehat G$ acts isometrically on $\widetilde X$ and this action centralizes the deck group $\pi=\pi_1(X)$.

First, we claim that the semisimple factor $S$ is actually trivial. Suppose that $S\neq\{1\}$ and choose a maximal torus $T\leq S$. Since $S$ is a nontrivial compact, connected Lie group, $\dim T>0$. Now, since $\Lambda$ is discrete in $\widehat G$, we have that $\Lambda\cap S$ is a closed discrete subgroup of the compact group $S$, hence finite.

The torus $T$ contains elements of order $2^r$ for arbitrarily large $r$.  Therefore there exists an element $s\in T$ of finite $2$-power order such that $s\notin\Lambda$. Thus $q(s)\neq e\in G$. Since the $G$-action on $X$ is effective, $q(s)$ acts non-trivially on $X$.  Hence its lift, namely the transformation induced by $s$ on
$\widetilde X$, is also nontrivial. The cyclic subgroup generated by this transformation is a nontrivial finite $2$-group of isometries of $\widetilde X$ centralizing $\pi$.  This contradicts Proposition \ref{prop:centralizer}. Therefore
$S=\{1\}$ and it follows that $\widehat G\cong\R^k$.

Once this claim is settled, we can observe that since $G$ is compact, $\Lambda\subset\R^k$ is a full lattice, and therefore $G\cong\R^k/\Lambda\cong\T^k$. Thus $G$ is a torus.

We now move on to prove the remaining thesis of the theorem, namely, that the morphism $\pi_1(G)\to Z(\pi)$ is injective. 
We now have
\[
    \Lambda= \Ker(\R^k\to G)\cong \pi_1(G) \cong \Z^k.
\]
As in Section \ref{s:lifts} there is a homomorphism $\rho :\Lambda\rightarrow Z(\pi)$. Suppose that $\rho$ is not injective and pick $0\neq\lambda\in\Ker\rho$. By definition of $\rho$, the element $\lambda\in\R^k$ acts trivially on $\widetilde X$. Consider the one-parameter subgroup $L:=\R\lambda\subset\R^k$. Since $\lambda$ acts trivially, the action of $L$ on $\widetilde X$ factors through $L/\Z\lambda\cong \mathbb{S}^1$. We claim that the resulting circle action on $\widetilde X$ is
nontrivial.

Indeed, if it were trivial, every element of $L$ would act trivially on
$\widetilde X$.  Since $\Lambda$ is discrete, there would exist $t\in\R$ such that $t\lambda\notin\Lambda$. Then $q(t\lambda)\neq e\in G$, while $q(t\lambda)$ would act trivially on $X$, contradicting the effectiveness of the $G$-action. It follows that the induced $S^1$-action on $\widetilde X$ is nontrivial.

Now, the kernel of this $\mathbb{S}^1$ action is a proper closed subgroup of $\mathbb{S}^1$ and every proper closed subgroup of $\mathbb{S}^1$ is finite. Hence, one can choose an element of $2$-power order in $\mathbb{S}^1$ outside this kernel. The cyclic group generated by its image is a nontrivial finite $2$-group acting  isometrically on $\widetilde X$. Since this circle action is induced by the lifted $\R^k$-action, it centralizes $\pi$. This again contradicts Proposition \ref{prop:centralizer}. This finally yields that $\Ker\rho=\{0\}$. 
\end{proof}

We now allow the acting compact Lie group $G$ to be disconnected.

\begin{proof}[Proof of Corollary \ref{cor:centerless}]
Let $K$ be a compact Lie group acting effectively and isometrically on $X$, and let $K^\circ$ denote its identity component. We observe that the restriction of the action $K^\circ\curvearrowright X$ is again effective. Then, by Theorem \ref{thm:main}, $K^\circ$ is a torus and $\pi_1(K^\circ)\hookrightarrow Z(\pi_1(X))$ is injective. 
Now, since $Z(\pi_1(X))=\{1\}$, we have a fortiori that $\pi_1(K^\circ)=0$ and a torus with trivial fundamental group is trivial, so $K^\circ=\{1\}$. Thus $K^\circ$ acts trivially on $X$. As the original $K$-action is effective, this forces $K^\circ=\{1\}$. Finally, a compact Lie group with trivial identity component is discrete, and a compact discrete group is finite.  Hence $K$ is finite.
\end{proof}

It is natural to wonder if the converse of Theorem \ref{thm:main} and Corollary \ref{cor:centerless} holds true. However the centerless hypothesis in Corollary \ref{cor:centerless} is sufficient, but not necessary, for the absence of symmetry. Indeed for every $n\geq2$ there exists a closed, connected, aspherical $n$-dimensional Alexandrov space $X$ such that $Z(\pi_1(X))\cong\Z^n$ and $\Isom(X)=\{1\}$ since for every compact smooth manifold of dimension at least two, metrics with trivial isometry group are generic (see, for example, the classical genericity theorem of Ebin \cite{Ebin}). The preceding examples only disprove the converse in the \emph{isometric} category.  Their underlying tori admit many
continuous torus actions.

Conner and Raymond asked whether nontrivial center should force a
circle action in the topological manifold setting.  Much later,
Cappell, Weinberger and Yan constructed closed aspherical manifolds with nontrivial center but no corresponding circle action, thereby showing that the topological converse also fails in general
\cite{CappellWeinbergerYan2013}.

% ============================================================
% Bibliography
% ============================================================
\bibliographystyle{siam}
\bibliography{references}

\end{document}